\documentclass[letterpaper, 10pt, conference]{ieeeconf}  

\IEEEoverridecommandlockouts                              
\usepackage{amsmath,amsbsy,amsgen,amscd,amssymb,amsfonts,stmaryrd}
\usepackage{mathtools}
\usepackage{url}
\usepackage{overpic}
\usepackage{subcaption}
\usepackage{lipsum}
\usepackage{multirow}
\usepackage[noadjust]{cite} 

\usepackage[dvipsnames]{xcolor}
\usepackage{pgfplots}
\pgfplotsset{compat=1.15}
\usepgflibrary{fpu}
\usetikzlibrary{calc, colorbrewer, arrows}
\usepgfplotslibrary{fillbetween}

\DeclareMathOperator*{\argmin}{arg\,min}
\DeclareMathOperator{\sech}{sech}

\title{\bf Learning latent representations in high-dimensional state spaces using polynomial manifold constructions}

\author{Rudy Geelen, Laura Balzano and Karen Willcox
\thanks{}
\thanks{Rudy Geelen is with the Oden Institute for Computational Engineering and Sciences, University of Texas at Austin, Austin, TX 78712 USA
        {\tt\small rudy.geelen@austin.utexas.edu}}%
\thanks{Laura Balzano is with the Department of Electrical Engineering and Computer Science, University of Michigan Ann Arbor, Ann Arbor, MI 48109 USA
        {\tt\small girasole@umich.edu}}%
\thanks{Karen Willcox is with the Oden Institute for Computational Engineering and Sciences, University of Texas at Austin, Austin, TX 78712 USA
        {\tt\small kwillcox@oden.utexas.edu}}%
}

\begin{document}

\maketitle
\thispagestyle{empty}
\pagestyle{empty}


\begin{abstract}

We present a novel framework for learning cost-efficient latent representations in problems with high-dimensional state spaces through nonlinear dimension reduction. By enriching linear state approximations with low-order polynomial terms we account for key nonlinear interactions existing in the data thereby reducing the problem's intrinsic dimensionality. Two methods are introduced for learning the representation of such low-dimensional, polynomial manifolds for embedding the data. The manifold parametrization coefficients can be obtained by regression via either a proper orthogonal decomposition or an alternating minimization based approach. Our numerical results focus on the one-dimensional Korteweg-de Vries equation where accounting for nonlinear correlations in the data was found to lower the representation error by up to two orders of magnitude compared to linear dimension reduction techniques.

\end{abstract}

\section{Introduction}
\label{sec:introduction}

Dimension reduction rests on the fundamental assumption that samples from high-dimensional state spaces can be represented with a smaller number of variables (often called ``latent" variables) without a significant loss of information. An arsenal of mathematical procedures has been proposed in the literature to find a compact representation of the data. The identification of such intrinsic, low-dimensional structure in problems with high-dimensional state spaces sits at the heart of this paper. In particular, we explore the use of state approximations that embed low-order polynomial terms. These terms are introduced to account for nonlinear correlations in the data thereby enabling a greater reduction in the dimensionality for the problem at hand. Our proposed approach remains interpretable through its formulation in terms of a modal basis expansion. As the least-squares method underlies the construction of these approximations, the framework can also be deployed effectively for large-scale, real-world datasets.

We start by summarizing the main linear technique for dimensionality reduction, namely, the principal component analysis (PCA) \cite{doi:10.1080/14786440109462720, hotelling1933analysis, jolliffe2005statistics}. The PCA performs a linear mapping of the data to a lower-dimensional space by maximizing its variance in the low-dimensional representation \cite{jolliffe2005statistics}. The breadth of applications for PCA, and its numerous adaptations, is largely due to its simplicity and ease of use. It should be noted, however, that PCA is sometimes also referred to as proper orthogonal decomposition (POD) \cite{lumley1967structures, sirovich1987turbulence, holmes1996turbulence}. PCA and POD are closely related and in some cases algorithmically equivalent, but have arisen from different communities. When it comes to dimension reduction in problems with high-dimensional state spaces associated with simulation data, POD is the method of choice. Both techniques have at their core the singular value decomposition (SVD) to identify low-rank structure. The methodology in this paper is presented from the perspective of POD, but everything we propose is equally applicable to PCA.

POD analysis dictates that the coefficients of the basis functions are linearly uncorrelated. Consequently, POD provides an optimal embedding of the original, high-dimensional data only into a linear subspace. It is unable to deal with nonlinear correlations existing in the data. The ability to define an adequate nonlinear embedding in high-dimensional state spaces is thus of great importance in dimension reduction problems. Many alternatives to POD have been proposed for identifying and unraveling such nonlinear correlations. These methods include, for instance, kernel PCA \cite{scholkopf2005kernel}, Isomap \cite{doi:10.1126/science.290.5500.2319}, locally linear embedding \cite{doi:10.1126/science.290.5500.2323}, multidimensional scaling (MDS) \cite{borg2005modern}, different variants of autoencoders, and other architectures that rely on neural networks as recently reviewed in, e.g., \cite{6472238}. Despite their successful deployment in many applications where linear dimensionality reduction techniques were found to fall short, machine learning methods are often criticized for their limited interpretability and susceptibility to overfitting. This can hinder, for instance, the embedding of physical constraints, such as boundary conditions, and the preservation of other data features and properties \cite{SWISCHUK2019704}. 

In this paper we introduce a new procedure for capturing key nonlinear interactions of the data by taking into account the existence of a low-dimensional, nonlinear manifold. The proposed strategy remains interpretable in that the nonlinear state approximations are given in terms of a modal basis expansion, as is the case with POD. The presented approach is also readily applicable in the dimension reduction of large-scale, real-world problems as its underlying optimization process is fully data-driven and can be carried out using standard least-squares solvers. The paper is organized as follows. First, we briefly review some known results for the POD as it applies to data compression problems. This is followed by a description of two instances of a nonlinear dimension reduction technique based on polynomial manifold constructions. The effectiveness of the proposed framework is then demonstrated in a nonlinear problem setting.

\section{Problem formulation}
\label{sec:pod}

The POD provides a systematic manner for the identification of a set orthonormal basis functions that can be used for approximating solutions in problems with high-dimensional state spaces. The POD basis is computed as the set of left singular vectors of a given data matrix $\mathbf{S} \in \mathbb{R}^{n \times k}$. The POD then seeks to produce rank-$r$ approximations of the state variable expressed through a linear combination of POD modes as follows:
\begin{equation}
    \mathbf{s}(t) \approx \boldsymbol{\Phi} \widehat{\mathbf{s}}(t),
    \label{eq:pod_approximation}
\end{equation}
where $\mathbf{s}(t) \in \mathbb{R}^n$ is the high-dimensional system state, $t$ is some parameter on which the state depends, $\boldsymbol{\Phi} \in \mathbb{R}^{n \times r}$ is a basis matrix (the POD basis) containing as columns the left singular vectors of $\mathbf{S}$ corresponding to the $r$ largest singular values. The column space of $\boldsymbol{\Phi}$ then defines an $r$-dimensional subspace of the full state space $\mathbb{R}^n$ with, generally speaking, $r \ll n$. The vector $\widehat{\mathbf{s}}(t) \in \mathbb{R}^r$ contains the reduced state (also POD) coordinates.

Provided that the total number of training samples is smaller than the dimension of the high-dimensional system states, that is $k < n$, the singular values of $\mathbf{S}$ are denoted $\sigma_1 \geq \dots \sigma_{k} \geq 0$. The POD basis minimizes the least-squares error between the original data and their representation in the reduced space (for a fixed basis size). This error is equal to the sum of the squares of the singular values corresponding to those left singular vectors not included in the basis
\begin{equation}
    \sum_{i=r+1}^{k} \sigma_i^2,
\end{equation}
where $\sigma_i$ is the $i$th singular value of $\mathbf{S}$. 

From a statistics viewpoint, the orthogonal projection underlying the POD provides linearly uncorrelated features. Despite its optimality properties, it cannot unravel nonlinear correlations in the data. Accounting for such nonlinear correlations is a crucial factor in further reducing the dimensionality of the problem. We now explore the construction of nonlinear state approximations that do not change the reduced problem dimension---in comparison to \eqref{eq:pod_approximation}---or call for additional degrees of freedom.

In the field of machine learning, the numerous tools aimed at identifying dominant nonlinear correlations are often referred to collectively as manifold learning or representation learning \cite{6472238}. We here focus on nonlinear transformations that can be expressed through a basis expansion. Of particular interest are nonlinear state approximations of polynomial form:
\begin{equation}
\mathbf{s}(t) \approx \underbrace{\mathbf{V} \widehat{\mathbf{s}}(t)}_\text{linear} + \underbrace{\mathbf{Z} \mathbf{g}(\widehat{\mathbf{s}}(t))}_\text{nonlinear},
\label{eq:nonlinear_approx}
\end{equation}
where
\begin{equation}
    \mathbf{g}(\widehat{\mathbf{s}}(t)) = \left[ \widehat{\mathbf{s}}^2(t), \widehat{\mathbf{s}}^3(t), \dots, \widehat{\mathbf{s}}^p(t) \right]^\top,
    \label{eq:polynomial}
\end{equation}
in which $\mathbf{V}  = [\mathbf{v}_1 \,| \, \dotsc \,|\, \mathbf{v}_r] \in \mathbb{R}^{n \times r}$ is the linear basis matrix that has $\mathbf{v}_j \in \mathbb{R}^n$, $j=1,\ldots,r$, as its $j$th column. The column space of $\mathbf{V}$ spans an $r$-dimensional linear subspace of the full state space $\mathbb{R}^n$. The matrix $\mathbf{Z} \in \mathbb{R}^{n \times (p-1)r}$ is a matrix whose columns span a different subspace from the one spanned by the columns of $\mathbf{V}$. The vector $\widehat{\mathbf{s}} \in \mathbb{R}^r$ is the reduced state vector, representing the coefficients of the basis expansion. The vector $\mathbf{g}(\widehat{\mathbf{s}}(t)) \in \mathbb{R}^{(p-1)r}$ denotes a set of polynomial terms up to degree $p \geq 2$ associated with the nonlinear aspect of the solution-manifold. Throughout this paper, the exponent used in \eqref{eq:polynomial} raises each element of a given vector to the corresponding power. For example, if $\widehat{\mathbf{s}}=[\widehat{s}_1,\widehat{s}_2]^\top$, then $\widehat{\mathbf{s}}^p = [\widehat{s}_1^p,\widehat{s}_2^p]^\top$. The use of polynomial constructions is not new in representation learning problems. However, the precise manner in which the polynomial approximations are built and deployed leads to vastly different approaches \cite{6220279, BARNETT2022111348, GEELEN2023115717, barnett2023mitigating, https://doi.org/10.48550/arxiv.2302.02036}.

If matrix $\mathbf{V}$ in \eqref{eq:nonlinear_approx} is chosen to be the POD basis matrix, $\mathbf{Z}$ can be inferred from the high-dimensional dataset by means of linear least-squares regression \cite{GEELEN2023115717}. A different approach is to consider a parametrization of the form
\begin{equation}
    \mathbf{Z} = \overline{\mathbf{V}} \boldsymbol{\Xi} \in \mathbb{R}^{n \times (p-1)r},
    \label{eq:parametrization}
\end{equation}
where the matrix $\overline{\mathbf{V}} = [\overline{\mathbf{v}}_1 \,| \, \dotsc \,|\, \overline{\mathbf{v}}_q] \in \mathbb{R}^{n \times q}$ represents a basis matrix, and $\boldsymbol{\Xi} \in \mathbb{R}^{q \times (p-1)r}$ is the corresponding coefficient matrix that weighs each of the polynomial terms contained in $\mathbf{g}(\widehat{\mathbf{s}}(t))$.

For optimization purposes, see Section \ref{sec:representation_learning}, we choose the columns of $\mathbf{V}$ and $\overline{\mathbf{V}}$ to form orthonormal sets of dimension $r$ and $q$, respectively. Additionally, the columns of $\overline{\mathbf{V}}$ are chosen to be orthogonal to those of $\mathbf{V}$. Equivalently we can write 
\begin{equation}
\lbrack \mathbf{V}, \overline{\mathbf{V}} \rbrack^\top \lbrack \mathbf{V}, \overline{\mathbf{V}} \rbrack = \mathbf{I}_{r+q},
\label{eq:orthogonality}
\end{equation}
where $\mathbf{I}_a$ denotes an identity matrix of dimension $a$. It is noted that constraint set \eqref{eq:orthogonality} forms a smooth submanifold of $\mathbb{R}^{n \times (r+q)}$ known as the Stiefel manifold.

To construct nonlinear approximations of the form \eqref{eq:nonlinear_approx} in an optimal fashion, the representation learning problem is posed in terms of the constrained optimization problem 
\begin{equation}
\begin{aligned}
\min_{\mathbf{V}, \overline{\mathbf{V}}, \boldsymbol{\Xi}, \widehat{\mathbf{S}}} & J(\mathbf{V}, \overline{\mathbf{V}}, \boldsymbol{\Xi}, \widehat{\mathbf{S}}) \\
\text{such that } & \lbrack \mathbf{V}, \overline{\mathbf{V}} \rbrack ^\top \lbrack \mathbf{V}, \overline{\mathbf{V}} \rbrack = \mathbf{I}_{r+q},
    \label{eq:optim_problem}
\end{aligned}
\end{equation}
with objective function
\begin{equation}
    J = \dfrac{1}{2} \sum_{j=1}^k \left\| \mathbf{s}_j - \lbrack \mathbf{V}, \overline{\mathbf{V}} \rbrack \begin{bmatrix} \widehat{\mathbf{s}}_j \\ \boldsymbol{\Xi} \mathbf{g}( \widehat{\mathbf{s}}_j ) \end{bmatrix} \right\|_2^2,
\end{equation}
where $ \widehat{\mathbf{S}} := (   \widehat{\mathbf{s}}_1 , \widehat{\mathbf{s}}_2 , \dots, \widehat{\mathbf{s}}_k ) \in \mathbb{R}^{r \times k} $ is the reduced-state representation of the given system states $\mathbf{s}_j$ for $j=1,\dots,k$. Solving \eqref{eq:optim_problem} means that 
\eqref{eq:nonlinear_approx} holds, in the least-squares sense, at the optimum of the objective of \eqref{eq:optim_problem}. 

The use of Frobenius norm regularization was advocated in \cite{GEELEN2023115717} to avoid the overfitting of the nonlinear approximations to the training data. Here we choose to penalise only the entries of the coefficient matrix $\boldsymbol{\Xi}$, leading to the modified optimization problem 
\begin{equation}
\begin{aligned}
    \min_{\mathbf{V}, \overline{\mathbf{V}}, \boldsymbol{\Xi}, \widehat{\mathbf{S}}} & \left( J(\mathbf{V}, \overline{\mathbf{V}}, \boldsymbol{\Xi}, \widehat{\mathbf{S}}) + \dfrac{\gamma }{2} \left\| \boldsymbol{\Xi} \right\|_F^2 \right) \\
    \text{such that } & \lbrack \mathbf{V}, \overline{\mathbf{V}} \rbrack^\top \lbrack \mathbf{V}, \overline{\mathbf{V}} \rbrack = \mathbf{I}_{r+q},
    \label{eq:optim_problem_reg}
\end{aligned}
\end{equation}
which is the problem of interest in this paper. Particular choices for each of the terms $\{\mathbf{V}, \overline{\mathbf{V}}, \boldsymbol{\Xi} \}$ and the corresponding latent space representation $\widehat{\mathbf{s}}(t)$, and how these may be inferred from the data, will be presented next.

\section{Learning nonlinear representations}
\label{sec:representation_learning}

\subsection{Approach 1 -- POD-based modeling}
\label{subsec:nonlinear_pod}

Our first option views the learning of nonlinear state approximations through the lens of the POD. In the following we demonstrate in step-by-step fashion how conventional POD approximations can be enriched to produce approximations of the form \eqref{eq:nonlinear_approx}.

\subsubsection*{Step 1: Compute basis}

\begin{figure}[tbp]
\centering
 \includegraphics[width=.95\linewidth]{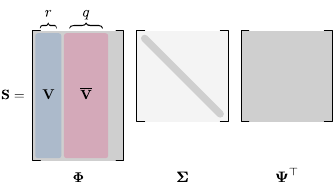}
 \caption{The left singular vectors used in the state approximations for proposed POD-based polynomial manifold approach following a singular value decomposition of the data matrix $\mathbf{S}=\boldsymbol{\Phi}\boldsymbol{\Sigma}\boldsymbol{\Psi}^\top$.}
 \label{eq:pod_approaches} 
\end{figure}

In this proposed approach we start by populating the basis matrix $\mathbf{V}$ with with the left singular vectors of $\mathbf{S}$ corresponding to its $r$ largest singular values. We then choose the basis vectors corresponding to the next consecutive $q$ largest singular values as the columns of $\overline{\mathbf{V}}$. A similar idea was also proposed recently in \cite{cohen:hal-04031976}. Figure \ref{eq:pod_approaches} illustrates the left singular vectors used in conventional POD and the proposed POD-based polynomial manifold approach.
Note that this choice for the pair of basis matrices \{$\mathbf{V}, \overline{\mathbf{V}}$\} automatically satisfies constraint \eqref{eq:orthogonality} through the orthogonality properties of the SVD. 

\subsubsection*{Step 2: Low-dimensional representation}

Finding the latent space representation for this approach simply amounts to computing the representation of each data sample $j$ in the POD coordinates as $\widehat{\mathbf{s}}_j = \mathbf{V}^\top \mathbf{s}_j$ for $j=1,\dots,k$.

\subsubsection*{Step 3: Linear least-squares regression}

Given the matrices $\{ \mathbf{V}, \overline{\mathbf{V}} \}$ and the latent space representation $\widehat{\mathbf{S}}$, this leaves the coefficient matrix $\boldsymbol{\Xi}$ to be inferred from the training data. It is then trivial to show that \eqref{eq:optim_problem_reg} reduces to an unconstrained linear optimization problem:
\begin{equation}
    \min_{\boldsymbol{\Xi}} \left( J(\mathbf{V}, \overline{\mathbf{V}}, \boldsymbol{\Xi}, \widehat{\mathbf{S}}) + \dfrac{\gamma }{2} \left\| \boldsymbol{\Xi} \right\|_F^2 \right),
    \label{eq:solving_for_psi}
\end{equation}
whose exact solution can be computed using the normal equations as
\begin{equation}
\boldsymbol{\Xi} = \overline{\mathbf{V}}^\top (\mathbf{I}_n-\mathbf{V}\mathbf{V}^\top)\mathbf{S} \mathbf{W}^\top ( \mathbf{W} \mathbf{W}^\top + \gamma \mathbf{I}_{(p-1)r})^{-1},
\label{eq:psi}
\end{equation}
with $\boldsymbol{\Xi} \in \mathbb{R}^{q \times (p-1)r}$, where we introduced the data matrix $\mathbf{W}$ as
\begin{equation}
    \mathbf{W} := 
    \begin{bmatrix}
        | & | & & | \\
        \mathbf{g}(\widehat{\mathbf{s}}_1) & \mathbf{g}(\widehat{\mathbf{s}}_2) & \dots & \mathbf{g}(\widehat{\mathbf{s}}_k) \\
        | & | & & |
    \end{bmatrix} \in \mathbb{R}^{(p-1)r \times k}.
\end{equation}
The vector accounting for the nonlinear interactions between the reduced-state coefficients of the $j$th data sample, $\mathbf{g}(\widehat{\mathbf{s}}_j)$, contains low-order polynomial nonlinear terms (up to degree $p$) and is given in \eqref{eq:polynomial}.

\subsubsection*{Contributions \& comparison with previous work}

The proposed POD-based polynomial manifold formulation described above produces only an approximate minimizer to optimization problem \eqref{eq:optim_problem}, as explained in \cite[Remark 1]{GEELEN2023115717}. Nonetheless, the formulation is particularly appealing from a computational perspective: it permits us to embed readily available POD modes explicitly into state approximation \eqref{eq:nonlinear_approx} without modifying the dimensionality of the reduced problem, an important contribution of this work. 

Despite also having its roots in polynomial regression, the POD-based manifold approach presented in \cite{GEELEN2023115717} is fundamentally different from the one described here in a number of key aspects. A parametrization of the form \eqref{eq:parametrization} permits us to explicitly leverage POD modes when approximating system states, rather than needing to compute the columns of $\mathbf{Z}$ from the data. As a result, the number of coefficients to be inferred in the proposed approach no longer scales with the full-state dimension, $n$, thereby reducing the computational and storage burden in constructing such approximations. Furthermore, the nonlinear part of the approximation is expressed through low-order polynomial terms instead of Kronecker quadratic product form. Accounting for higher-order nonlinear correlations existing between the various POD mode amplitudes will promote an improved compression of the data.

\subsection{Approach 2 -- The alternating minimization approach} 
\label{subsec:alternate_minimization}

In this section we take a different approach from above in that we propose an iterative algorithm for solving problems of the form \eqref{eq:optim_problem_reg} using alternating minimization. Classical alternating minimization (or projection) algorithms have been successfully applied in solving optimization problems over two or more variables or equivalently of finding a point in the intersection of two sets. These methods have a long history in many areas such as signal processing, control, and finance \cite{4557457}. Although one sacrifices a certain degree of interpretability by solving the representation learning problem through numerical optimization procedures, each step of the proposed method is intuitive yet computationally tractable and thus scalable. Our approach is detailed in individual subsections below in step-by-step fashion. We also provide an adequate choice of the initialization as well as a relevant stopping criterion in the alternating minimization process.

\subsubsection*{Step 1: The Procrustes problem}

We start with initial guesses for the representation of the training data in the low-dimensional coordinate system at each sample $\widehat{\mathbf{s}}_j$, $j=1,\dots,k$, and the coefficient matrix $\boldsymbol{\Xi}$. The matrix of the concatenated bases $\boldsymbol{\Omega} := [\mathbf{V} , \overline{\mathbf{V}} ]$ can then be inferred from the available data. Because the regularizer in \eqref{eq:optim_problem_reg} depends neither on $\mathbf{V}$ nor on $\overline{\mathbf{V}}$, an equivalent optimization problem may be considered without regularization. In the Frobenius norm \eqref{eq:optim_problem_reg} then simplifies to
\begin{equation}
\begin{aligned}
    \boldsymbol{\Omega} = \argmin_{\boldsymbol{\Omega} \in \mathbb{R}^{n \times (r+q)}} & \dfrac{1}{2} \left\| \mathbf{S} - \boldsymbol{\Omega} \begin{bmatrix} \widehat{\mathbf{S}} \\ \boldsymbol{\Xi} \mathbf{W}\end{bmatrix} \right\|_F^2 \\ \text{such that } & \boldsymbol{\Omega}^\top \boldsymbol{\Omega} = \mathbf{I}_{r+q}.
    \label{eq:procrustes}
    \end{aligned}
\end{equation}
Matrix approximation problem \eqref{eq:procrustes} is better known in the statistical literature as the orthogonal Procrustes problem \cite{gower2004procrustes}. It aims to learn the best possible unitary (orthogonal) transformation, that is a rotation or reflection, that relates two given matrices. The solution to \eqref{eq:procrustes} was derived by Sch\"{o}nemann \cite{schonemann1966generalized} as 
\begin{equation}
\boldsymbol{\Omega} = \mathbf{U}_{P}\mathbf{V}_{P}^\top,  
\end{equation}
 where $\mathbf{U}_{P}\boldsymbol{\Sigma}_{P}\mathbf{V}_P^\top = \mathbf{S}[ \widehat{\mathbf{S}}^\top, (\boldsymbol{\Xi} \mathbf{W})^\top ]$ is an SVD, and is unique if $\mathbf{S}[\widehat{\mathbf{S}}^\top, (\boldsymbol{\Xi} \mathbf{W})^\top]$ is full rank.

\subsubsection*{Step 2: Linear least-squares regression}

Given an estimate for the basis matrices, $\{ \mathbf{V}, \overline{\mathbf{V}} \}$, and the same low-dimensional coordinate system used in Step 1, we now solve for the coefficient matrix $\boldsymbol{\Xi}$. In practice this amounts to solving \eqref{eq:solving_for_psi}, with the main difference being one of interpretation: in Approach 1 the basis matrices $\{ \mathbf{V}, \overline{\mathbf{V}} \}$ are the left singular vectors of the high-dimensional training data whereas in the alternating minimization approach this is no longer the case as one has to solve an orthogonal Procrustes problem to find a set of orthonormal basis matrices.

\subsubsection*{Step 3: Unconstrained nonlinear least-squares}

In the final step of the minimization process we compute the representation of the high-dimensional data on the polynomial manifold. For this we have new estimates for $\mathbf{V}, \overline{\mathbf{V}}$ and $\boldsymbol{\Xi}$ at hand, meaning that \eqref{eq:optim_problem_reg} simplifies to
\begin{equation}
\begin{aligned}
    \widehat{\mathbf{S}} = & \argmin_{ \widehat{\mathbf{S}} \in \mathbb{R}^{r \times k} } \dfrac{1}{2} \sum_{j=1}^k  \left\| \mathbf{s}_j - \lbrack \mathbf{V}, \overline{\mathbf{V}} \rbrack \begin{bmatrix} \widehat{\mathbf{s}}_j \\ \boldsymbol{\Xi} \mathbf{g}( \widehat{\mathbf{s}}_j ) \end{bmatrix} \right\|_2^2,
    \label{eq:opt_low_dim_coord_system}
\end{aligned}
\end{equation}
which is an unconstrained optimization problem that can thus be treated efficiently using Levenberg-Marquardt and trust-region-reflective methods available in many off-the-shelf solvers.

\subsubsection*{Initial guess and termination criterion}
\label{subsec:termination}
The iterative process must be initiated with guesses for the representation of the data in the low-dimensional coordinate system, $\widehat{\mathbf{s}}_j$ for $j=1,\dots,k$, as well as the coefficient matrix $\boldsymbol{\Xi}$. These are readily available in our POD-based polynomial manifold approach (Approach 1). In the proposed alternating minimization scheme we solve the problems in Steps 1 to 3 in a sequential fashion, constituting one iteration in the optimization procedure. Iterations are terminated when a global stopping criterion is satisfied. We choose to monitor the retained information content, see \cite{GEELEN2023115717} for a motivation, accounted for by approximation \eqref{eq:nonlinear_approx} across the iterations, denoted by the $\ell$ superscript, as
\begin{equation}
    e^{\ell} := \dfrac{\| \mathbf{V}^{\ell}\widehat{\mathbf{S}}^{\ell} + \overline{\mathbf{V}}^{\ell}\boldsymbol{\Xi}^{\ell}\mathbf{W}^{\ell} \|_F^2}{ \| \mathbf{S} \|_F^2 }.
    \label{eq:convergence}
\end{equation}
We then define convergence if $| e^{\ell+1}-e^{\ell} | \leq \texttt{TOL}$ with $\texttt{TOL}$ a predefined tolerance.

\section{Numerical results}
\label{sec:numerical_results}

\subsection{Error metric \& good practices}
\label{subsec:error_metrics}

In the numerical experiments conducted in this work we represent a given dataset $\mathbf{S}(\boldsymbol{\mu}_i)$ in the computed reduced-state coordinate system, after which we attempt to reconstruct that same dataset in the original, high-dimensional state space as $\widetilde{\mathbf{S}}(\boldsymbol{\mu}_i)$. The different parameter instances are given by $\boldsymbol{\mu}_i$ with $i=1,\dots,m$. Importantly, the data used for testing are unseen in learning process. This enables us to quantify the representation error of nonlinear state approximation \eqref{eq:nonlinear_approx} as follows:
\begin{equation}
    \dfrac{1}{m}\sum_{i=1}^m \dfrac{\| \mathbf{S}(\boldsymbol{\mu}_i)-\widetilde{\mathbf{S}}(\boldsymbol{\mu}_i) \|_F^2 }{ \| \mathbf{S}(\boldsymbol{\mu}_i) \|_F^2 }
    \label{eq:representation_error}
\end{equation}

The data matrices under consideration are centered by their column-averaged mean value. The regularization parameter $\gamma$ in \eqref{eq:optim_problem_reg} is found through a grid search and confirmed to provide an adequate amount of regularization irrespective of the polynomial degree of the state approximation or the number of basis vectors $r+q$. The decay of the singular values is used in making an informed choice for $q$. The tolerance for the alternating minimization algorithm, see Section \ref{subsec:termination}, is chosen to be $\texttt{TOL}=10^{-3}$.

\subsection{Korteweg-de Vries equation}
\label{subsec:kdv}

Our numerical experiments are concerned with traveling wave physics. Consider a single soliton propagating at constant speed over a domain with periodic boundary conditions. The evolution of the wave field is given by the one-dimensional Korteweg-de Vries equation
\begin{equation}
    \partial_t s = -\alpha s \partial_x s - \beta \partial_{xxx} s
\end{equation}
in the space-time domain $[-\pi,\pi]\times[0,0.1]$ with initial condition $s_0(x) = 1 + 24 \sech^2 \left( \sqrt{8}(x-\mu) \right)$ and with $\partial_x$ and $\partial_t$ denoting partial derivatives with respect to $x$ and $t$, respectively \cite{mendible2020dimensionality}. An equidistant grid was used consisting of 256 grid points. The state data are saved every $0.0002$ time units. We choose the model constants to be $\alpha=8$ and $\beta=1$ and generate five individual datasets at $\mu = [0,0.5,1,1.5,2]$, yielding a total of 2,500 samples to be used for learning the state representations.
\begin{figure}[tbp]
\centering \footnotesize
\includegraphics{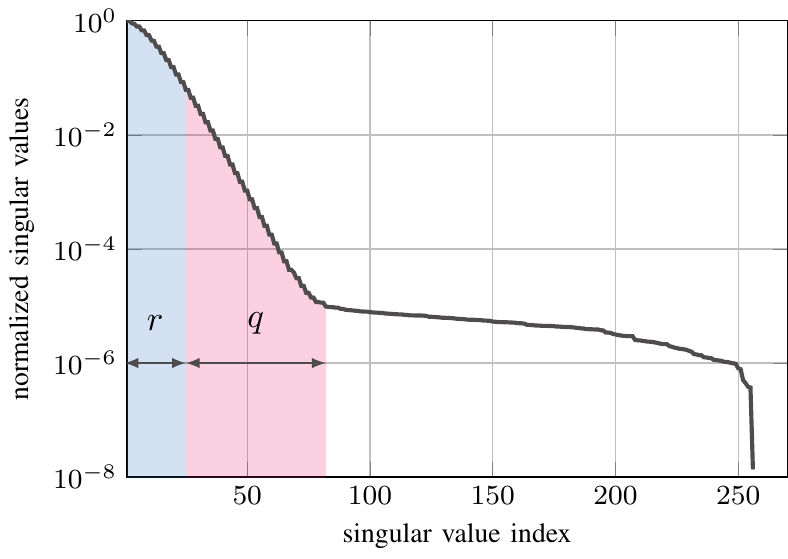}
\caption{Normalized singular values of the mean-subtracted data matrix for the Korteweg-de Vries problem. The blue and red areas denote the singular values whose corresponding left singular vectors are in $\mathbf{V}$ and $\overline{\mathbf{V}}$, respectively.}
\label{fig:svs_kdv}
\end{figure}

As evidenced by the decay of the singular values of the dataset, shown in Fig.\ \ref{fig:svs_kdv}, 82 POD basis vectors are necessary for approximating the training data with a projection error below $10^{-5}$ in the Euclidean norm. The decay is also found to reach a plateau around this value, at which point increasing the reduced basis dimension no longer meaningfully improves the quality of the state approximation. In the numerical experiments that follow we therefore maintain a constant $r+q=82$. The regularization parameter was chosen to be $\gamma=500$. 

\begin{figure}[tbp]
\centering \footnotesize
\begin{subfigure}[t]{.49\linewidth}
    \includegraphics[width=\linewidth]{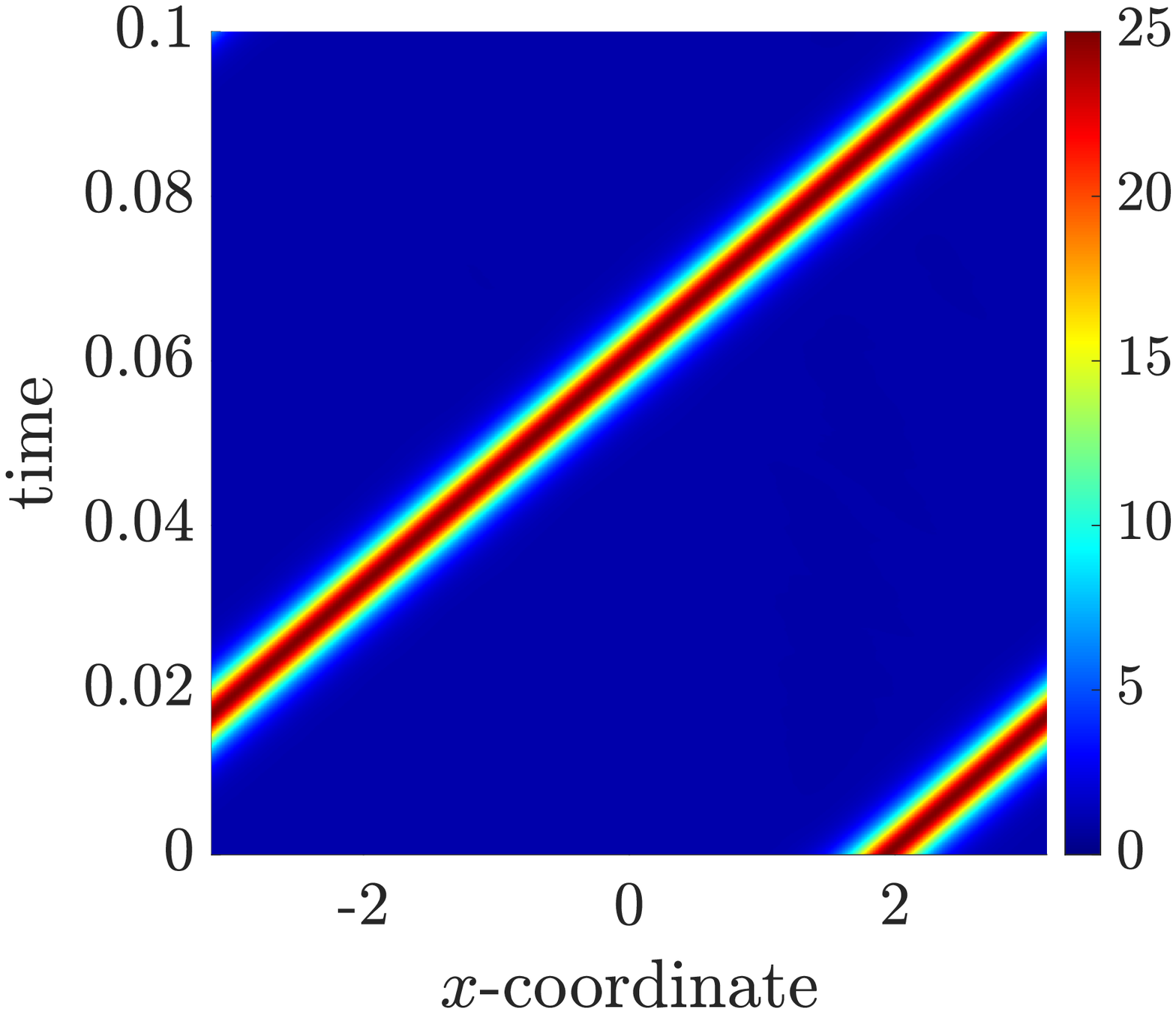}
    \caption{\footnotesize Reference solution}
\end{subfigure}
\begin{subfigure}[t]{.49\linewidth}
    \includegraphics[width=\linewidth]{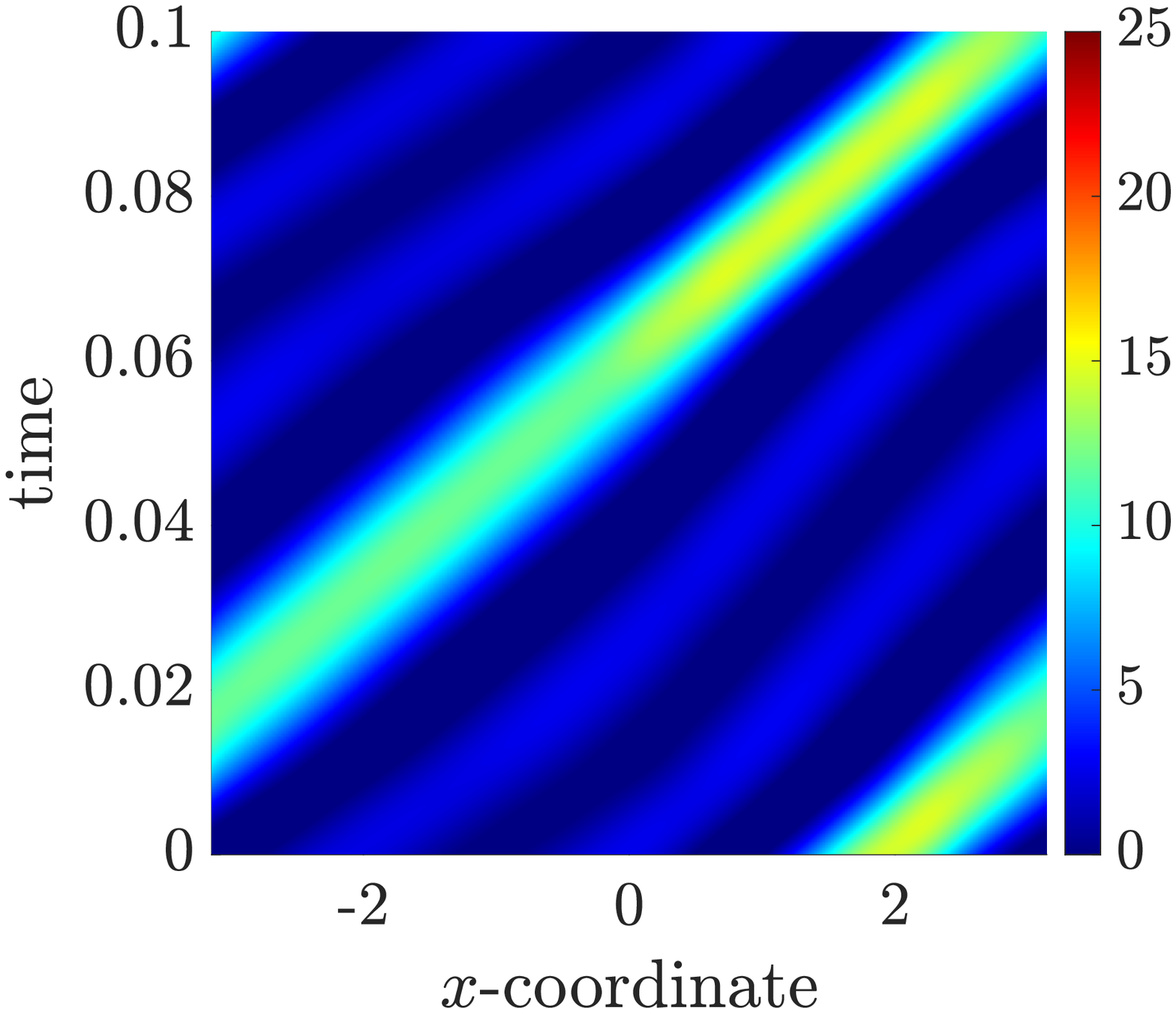}
    \caption{\footnotesize POD reconstruction with $r=6$}
\end{subfigure}\\[1em]
\begin{subfigure}[t]{.49\linewidth}
    \includegraphics[width=\linewidth]{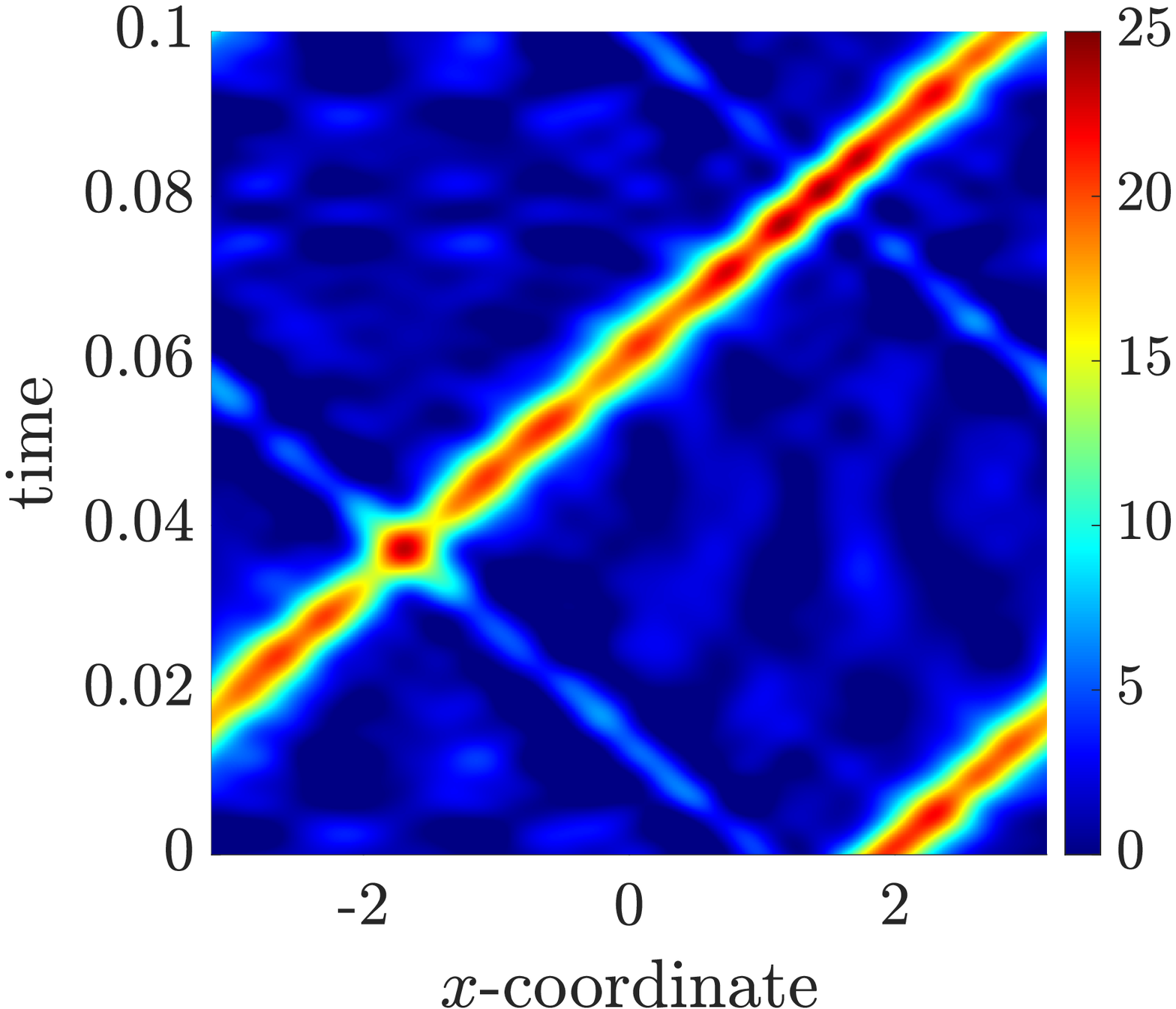}
    \caption{\footnotesize Manifold reconstruction (POD) with $r=6, p=4$}
    \label{fig:kdv1}
\end{subfigure}
\begin{subfigure}[t]{.49\linewidth}
    \includegraphics[width=\linewidth]{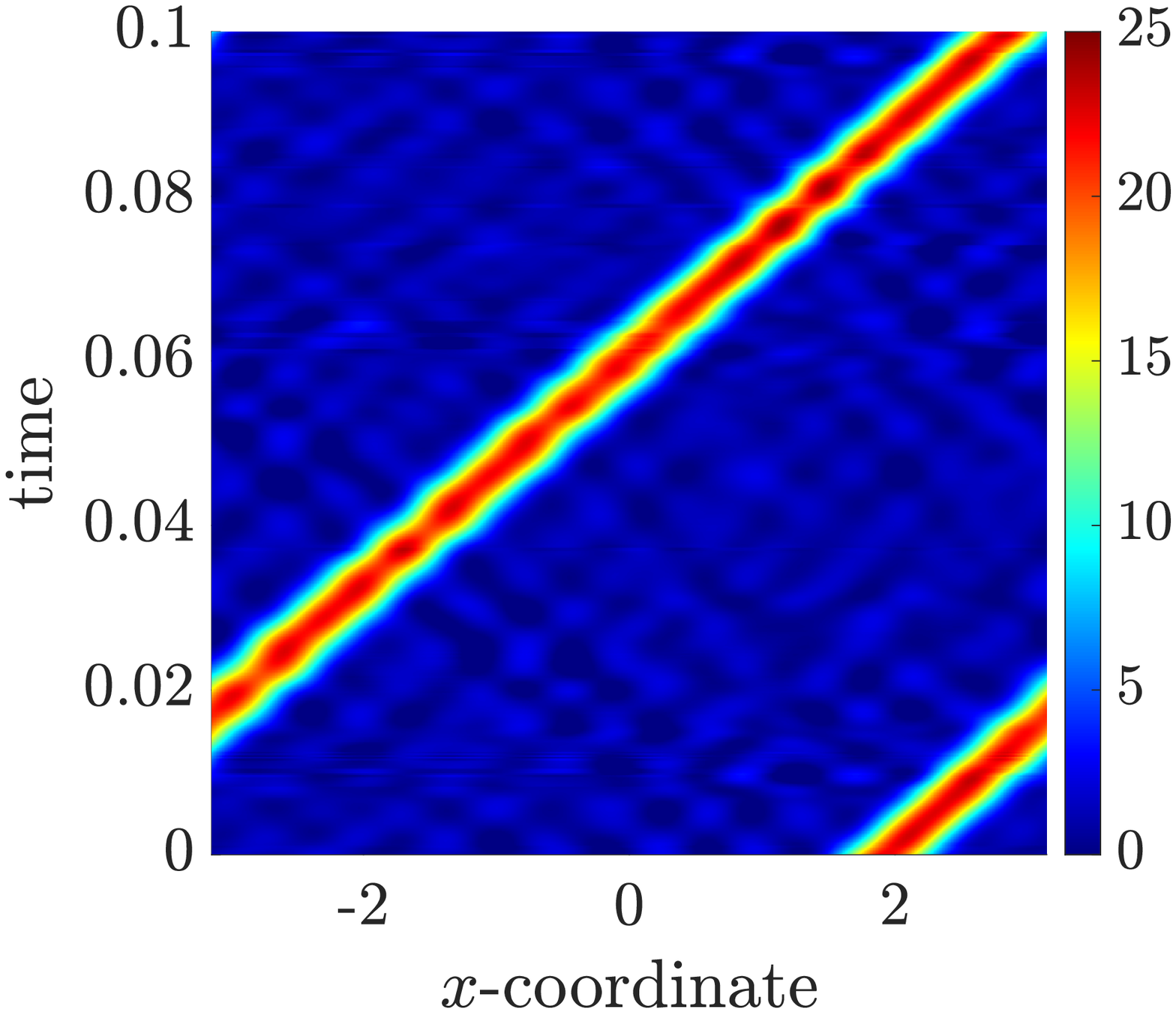}
    \caption{\footnotesize Manifold reconstruction (AM) with $r=6, p=4$}
    \label{fig:kdv2}
\end{subfigure}
\caption{Reconstruction of the high-dimensional system-states using the different state approximation at $\mu_j = 1.9298$ with reduced basis dimension $r=6$. The manifold models have polynomial embeddings of degree $p=4$.}
\label{fig:kdv_error_space_time}
\end{figure}

An additional ten datasets are generated, based on the initial conditions, by uniformly drawing $\mu$ in the range $[0,2]$. In this experiment we are looking to better understand the ability of the POD, Manifold-POD (Section \ref{subsec:nonlinear_pod}), and Manifold-AM (Section \ref{subsec:alternate_minimization}) formulations in building an $r$-dimensional reduced-space representation and, from that, reconstruct a set of high-dimensional system-states with a minimal loss of information. Fig.\ \ref{fig:kdv_error_space_time} gives an illustration of the performance of the different reduction methods at a reduced dimension of $r=6$.
The numerical experiment considers a polynomial embedding of degree $p=4$ for the nonlinear state approximations. The POD clearly struggles with capturing the space-time evolution of the Korteweg-de Vries equation. As it represents the optimal solution that can be approximated using a linear subspace spanned by only $r$ basis vectors, the reconstruction fails to capture the soliton's amplitude and shows Gibbs oscillations in its vicinity. The reconstructions from the polynomial manifold formulations on the other hand, shown in Figs.\ \ref{fig:kdv1} and \ref{fig:kdv2}, are clearly better suited for this by accounting for correlations of polynomial nature between the latent space coefficients. This ability is beneficial in reducing the problem's dimensionality.

\begin{figure}[tbp]
\footnotesize \centering
\includegraphics{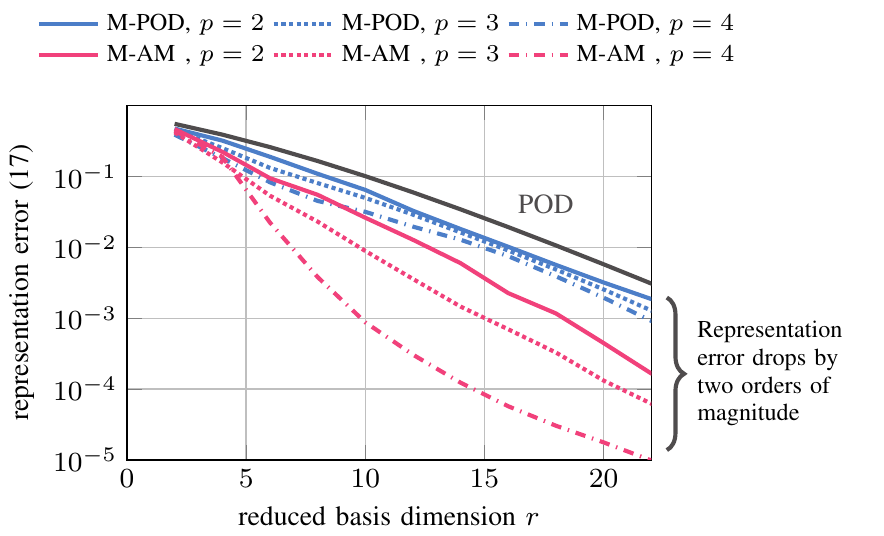}
\caption{Plot of the representation error, \eqref{eq:representation_error}, in function of the reduced basis dimension, $r$, for the standard POD approach (black curve), POD-based polynomial manifold approach (blue curves) and alternating minimization based polynomial manifold approach (red curves). For the manifold formulations we vary the degree $p$ of the polynomial embedding.}
\label{fig:kdv_error}
\end{figure}

The test errors in function of the reduced basis dimension are shown in Fig.\ \ref{fig:kdv_error}.
As expected, POD is the least able in accurately reconstructing the test data. When accounting for an additional $q$ POD modes through the addition of low-order polynomial terms to the state approximation, the representation error decreases. However, the proposed alternating minimization approach, in which we rely on the solution to a Procrustes-type problem for inferring a set of orthonormal basis vectors, is considerably more accurate than both aforementioned alternatives. Improvements up to two orders of magnitude in accuracy were obtained compared to POD. Both nonlinear manifold approaches were found to benefit from increasing the polynomial degree in the nonlinear part of the approximation \eqref{eq:nonlinear_approx}.

The number of iterations at which \eqref{eq:convergence} falls below the specified threshold for the alternating minimization procedure is listed in Table \ref{table_example}. While there is a quite a bit of variation at different values of the reduced basis dimension, $r$, and the degree of the polynomial embeddings, $p$, the number of iterations remains small for all numerical experiments reported here. This shows that learning nonlinear state approximations of the form \eqref{eq:nonlinear_approx} through minimization techniques is computationally tractable.

\begin{table}[tbp]
\caption{The number of alternating minimization cycles at which \eqref{eq:convergence} falls below the specified threshold at different reduced basis dimensions, $r$, and values for the polynomial embedding degree, $p$.}
\label{table_example}
\begin{center}
\begin{tabular}{|l|l|c|c|c|c|c|c|c|}
\hline
\multicolumn{2}{|l|}{Reduced basis dim.\ $r$} & $2$ & $4$ & $6$ & $8$ & $10$ & $12$ & $14$ \\ \hline
\multirow{3}{*}{\# iterations}  & $p=2$  & 4  & 8  & 11 & 20 & 9 & 10 & 5 \\
                                & $p=3$  & 6  & 13 & 18 & 15 & 6 & 5 & 4 \\
                                & $p=4$  & 19 & 29 & 9  &  6 & 4 & 3 & 3 \\
\hline
\end{tabular}
\end{center}
\end{table}

\section{Concluding remarks}
\label{sec:conclusions}

We present a general framework for building nonlinear state approximations that is particularly applicable for dimension reduction in problems with high-dimensional state spaces. By embedding low-order polynomial terms in the modal basis expansion we account for nonlinear correlations in the data through interactions of the reduced space coordinates. The proposed framework prompts two different representation learning approaches. The POD-based polynomial manifold approach is intuitive due its connection to conventional POD, yet enables a reduction of the dimensionality for the problem at hand. If one is willing to depart from POD, and its interpretable nature, the proposed alternating minimization based variant offers an even better compression performance by finding basis matrices with orthonormal columns through a Procrustes-type problem. Importantly, both learning approaches can effectively be deployed at scale using only least-squares techniques.

\section{Acknowledgements}

This work has been supported in part by the U.S.\ Department of Energy AEOLUS MMICC center under award DE-SC0019303, program manager W.\ Spotz, and by the AFOSR MURI on physics-based machine learning, award FA9550-21-1-0084, program manager F.\ Fahroo. Laura Balzano was supported by ARO YIP award W911NF1910027 and NSF CAREER award CCF-1845076. Karen Willcox would also like to thank Yvon Maday and Albert Cohen for several helpful conversations during the 2022 Leçons Jacques-Louis Lions in Paris, France.

\bibliographystyle{ieeetr}
\bibliography{main}

\end{document}